\documentclass[11pt]{amsart}
\usepackage{amssymb}
\usepackage{amsmath}
\usepackage{amsfonts}
\usepackage{color}
\newtheorem{theorem}{Theorem}
\theoremstyle{plain}
 
\newtheorem {definition}[theorem]{Definition}
\newtheorem {lem}{Lemma}
\newtheorem {prop}{Proposition}
\newtheorem {remark}{Remark}
\newtheorem {Rem}{Remark}
\newtheorem {theo}{Theorem}
\newtheorem {corol}{Corollary}
\newtheorem {example}{Example}
\begin{document}
\renewcommand{\bibname}{Bibliography}
\title{ On some properties of $n$-EP and $n$-hypo-EP operators. }
\author{R. Semmami}
\address{R. Semmami,  Laboratory of Mathematical Analysis and Applications, Faculty of Sciences, Mohammed V University in Rabat, Rabat, Morocco.}
\email{semmami.rachi@gmail.com}
\author{H. Ezzahraoui}
\address{H. Ezzahraoui, Laboratory of Mathematical Analysis and Applications, Faculty of Sciences, Mohammed V University in Rabat, Rabat, Morocco.}
\email{hamid.ezzahraoui@fsr.um5.ac.ma;h.ezzahraoui@um5r.ac.ma}
\author{E. H. Zerouali}
\address{Permanent address: Laboratory of Mathematical Analysis and Applications, Faculty of Sciences, Mohammed V University in Rabat, Rabat, Morocco.\\
Current address: Department of Mathematics, The University of Iowa, Iowa City, Iowa, U.S.A.}
\email{elhassan.zerouali@fsr.um5.ac.ma \& ezerouali@uiowa.edu}
\thanks{The last author was partially supported by the Arab Fund Foundation Fellowship Program.
The Distinguished Scholar Award - File  1026. He also acknowledges the mathematics department of the University of Iowa for kind hospitality during the preparation of this paper.}

\begin{abstract}
This article, is devoted to  $n$-EP and  $n$-hypo-EP operators. We  give some  characteristic properties of these two classes and   various links with other known classes in the literature, especially the classes of EP, SD, hypo-EP and $n$-normal operators.
\end{abstract}
\maketitle
\section{Introduction}
Let $\mathcal{H}$ be a Hilbert space and let us denote  $\mathcal{L(H)}$ for the space of
all bounded linear operators defined on   $\mathcal{H}$ and let $\mathcal{R(H)}\subset \mathcal{L(H)}$ be the set of  operators with closed range. For $T \in \mathcal{L(H)}$, we denote by $R(T)$, $N(T)$ and $T^{*}$ the range of $T$, the kernel of $T$ and the adjoint operator of $T$ respectively.

A closed subspace $E\subset\mathcal{H} $ is said to be \emph{invariant} for $T$ if $T(E) \subset E$ and it is called  \emph{  reducing  }
 for $T$ if it is invariant for both $T$ and $T^*$, or equivalently, if $T(E)\subset E$ and $T(E^{\perp})\subset E^\perp$, where  $E^\perp$ is the orthogonal complement of $E$.

  For $T,S\in \mathcal{L(H)}$, we write $[T,S]:=TS-ST$ for the commutator of $S$ and $T.$ The operator    $T$ is called normal if   $[T^{*},T]=0,$  is quasi-normal if $[T^*T,T]=0$ and  is hyponormal if $[T^{*},T]\ge 0.$

 An operator $T\in\mathcal{L}(\mathcal{H})$ is  a \emph{partial isometry} if $TT^*T=T$, equivalently $T^*T= P_{N(T)^\perp}$.   Partial isometries have been investigated in \cite{mbk3}. It is known   that $T$ is a partial isometry if and only if $T^*$ is a partial isometry

For   $T\in \mathcal{R(H)}$, the Moore-Penrose inverse $T^{\dagger}$ of $T$ is defined as the unique solution of the following four operator equations: \begin{equation}\label{morequ}
    TT^{\dagger}T=T,\;T^{\dagger}TT^{\dagger}=T^{\dagger},\;(T^{\dagger}T)^{*}=T^{\dagger}T\;\text{and }(TT^{\dagger})^{*}=TT^{\dagger}.
\end{equation}
The next properties of the Moore-Penrose inverse are immediate from Equations \ref{morequ}, and  are to be used in this paper.
 \begin{prop}\label{moreprop} If $T\in\mathcal{R(H)} $, then
 \begin{enumerate}
     \item   $TT^{\dag}=P_{R(T)}$, $T^{\dag}T= P_{R(T^{\dag})}$,
     \item  $R(T^{\dag})=R(T^*)$  and  $N(T^{\dag})=N(T^*)$,
     \item  $T^{*}TT^{\dag}=T^{\dag}TT^{*}=T^{*}$,
     \item  $T$ is a partial isometry if and only if $ T^{*}=T^{\dag}$,
     \item  $T^\dag = \lim\limits_{s\to 0}(T^*T+sI)^{-1}T^*.$
 \end{enumerate}
 \end{prop}
For  equivalent definitions and further properties of $T^{\dagger}$, we refer to \cite{Adi, sameer1, emz1,  katzand1}.

 An operator $T\in\mathcal{R(H)} $  is called an EP operator ($T\in \text{(EP)}$, for short) if $R(T)=R(T^{*}).$ This concept was introduced for matrices by Schwerdtfeger in \cite{schwerdtfeger} and has been studied in detail by several authors. See  \cite{baskett, katzand1},  for example.
It is easy to verify that
$$T\in \text{(EP)} \iff T^*\in \text{(EP)} \iff T^{\dagger}\in \text{(EP)} \iff [T^{\dagger}, T] = 0.$$
Recall that $T$ is said to be left invertible if $ST=I$ for some bounded operator $S$. It is known that $T$ is left invertible if and only if $T$ is one to one and has closed range. It follows then that a left invertible operator is an EP if and only if, it is invertible.

The Cauchy dual of a left invertible operator $T\in\mathcal{L(H)}$ is given by $\omega(T)=T(T^{*}T)^{-1}= (T^{\dag})^{*}$ and is an efficient tool in the study of structural properties of left invertible operators.  From the definition,  a left invertible operator $T$ is isometric,  if and only it $\omega(T)=T$.

The notion of Cauchy dual has been extended to    closed range operators in \cite{emz1}.  More precisely, the Cauchy dual of an operator $T\in\mathcal{R(H)}$ is defined as $\omega(T)=T(T^{*}T)^{\dagger}.$
It is easily seen that
$$ \omega(T)=  T^{*\dag}
           = T(T^{*}T+P_{N(T)})^{-1}
           =  (TT^{*}+P_{N(T^*)})^{-1}T,
$$
where   $P_E$ denotes the orthogonal projection on a given subspace $E$. In particular $T$ is a partial isometry if and only if $\omega(T)=T$.
  It is also easy to see that  $N(\omega(T))=N(T)$ and $R(\omega(T))=R(T)$. The reader is referred  to \cite{emz1} for more information.

In the line of hyponormality, M. Itoh introduced in \cite{itoh} the class of hypo-EP operators ((HEP) for short)   as follows. An operator $T\in \mathcal{R(H)}$ is an  hypo-EP operator if  $[T^{\dagger},T]\geq 0$   or equivalently $R(T)\subset R(T^{*})$. See \cite{ding1, vinoth} for detailed study and further results on HEP operators.

An operator $T\in\mathcal{R(H)}$ is  SD, the class of star-dagger operators, if  $[T^*,T^\dag]=0$, or equivalently $[T,\omega(T)]=0$. The special class of SD matrices appeared in \cite{hartwig} and has been widely studied in \cite{ferr,hartwig}.

Our purpose in this paper is to explore some relationships between EP operators, normal operators, SD operators  and HEP operators. In particular, we provide a characterization of such  operators and we give necessary and sufficient conditions  for HEP operators to be EP operators.

The concept of $n$-EP operators  is introduced in  \cite{SLN} for matrices and extended in \cite{wangchun}  to bounded operators with closed range where several properties of $n$-EP operators  are provided. For further properties of $n$-EP operators, we refer to \cite{emz1,menkad}. We introduce in this paper a new class of $n$-hypo-EP operators  extending  (EP) and (HEP) classes. We give some additional properties of $n$-EP operators and we study the class of HEP operators.
\section{Some proprieties of $n$-EP operators}
\subsection{EP and SD operators}
Before exposing additional results concerning the  classes (SD) and (EP), we give some examples of EP and SD operators.
\begin{example}\cite[Example 2.3]{emz1}
  Let $T=x\otimes y$ be a one rank  operator on $\mathcal{H}$. Since $R(T)=\mathbb{C}x$ is closed and $R(T^*)=\mathbb{C}y$, it follows that,  $T$ is an EP operator if and only if $y=\alpha x$ for some nonzero $\alpha\in \mathbb{C}$.

  On the other hand, since the Cauchy dual of $T$ is   $$\omega(T)=\frac{1}{\|x\|^{2}\|y\|^{2}}x\otimes y=\frac{1}{\|x\|^{2}\|y\|^{2}}T,$$
  we derive that $T$ is always an SD operator.
\end{example}
\begin{example}
  Let $T$ be a bounded orthogonal projection on a Hilbert space $\mathcal{H}$ (i.e. $T^{2}=T=T^{*}$). We have $\omega(T)=T= T^\dag$, thus $T\in\text{(EP)} \cap\text{(SD)} $.
\end{example}
\begin{example}
 Let $\mathcal{H}$ be a Hilbert space endowed with some orthonormal basis $\{e_{n}\}_{n\geq0}$. The unilateral weighted shift associated with some non negative sequence $(\alpha_n)_{n\ge 0}$ is the linear operator defined on the
basis of $\mathcal{H}$ by $S_{\alpha}e_{n} = \alpha_{n} e_{n+1}$ for all $n\geq0.$  The adjoint operator of $S_{\alpha}$ is given by $S_{\alpha}^{*}e_{n} = \alpha_{n-1} e_{n-1}$ if $n\geq 1$ and $S^{*}_{\alpha} e_{0}=0.$ It is known that $S_\alpha$ has closed range if and only if
$\lim\limits_{n\rightarrow+\infty}(\inf\limits_{k\ge 0}\alpha_k\cdots\alpha_{k+n-1})^{\frac{1}{n}}>0.$
In this case  $R(S^*_\alpha)= \mathcal{H}\ne \{e_{0}\}^\perp =R(S_\alpha)$ and then  $S_\alpha$ is never EP.

On the other hand, direct computations show that  $\omega(S_{\alpha})$ is the unilateral weighted shift defined by
$$  \omega(S_{\alpha})e_{n}=\frac{1}{\alpha_{n}}e_{n+1} \mbox{ for every } n.$$

We obtain,  $S_{\alpha}$ is SD, if and only if   $\omega(S_\alpha )S_\alpha= S_\alpha\omega(S_\alpha )$ and equivalently
$\frac{\alpha_{n+1}}{\alpha_{n}}=\frac{\alpha_{n}}{\alpha_{n+1}}$ for every $n\ge 0$. It follows that $S_{\alpha}$ is SD if and only if  $\alpha_n=\alpha_0$ for every $n\ge 1$.
\end{example}
Recall that  $S$ and   $T$ are said to be unitarily equivalent if there exists a unitary operator $U$ such that $S=U^*TU.$ The class $(SD)$ is closed under unitary equivalence.
\begin{theo}
 Let $S$ and   $T$   be unitarily equivalent. Then  $$T\in\text{(SD)}  \iff S\in \text{(SD)} $$
\end{theo}
\begin{proof}
It is clear that $R(S)$ is closed, since $R(T)$ is closed. From \cite[Proposition 2.2]{emz1}, we obtain $\omega(S)= ((U^*TU)^*)^\dag = U^*\omega(T)U$. Therefore,
\begin{eqnarray*}
\omega(S)S &=&U^*\omega(T)U.U^*TU\\
           &=&U^*\omega(T)TU\\
           &=&U^*T\omega(T)U\\
           &=&U^*TU.U^*\omega(T)U\\
           &=&S\omega(S).
\end{eqnarray*}

\end{proof}
Clearly any invertible operator is EP, while an invertible operator is  SD if and only if it is  normal. Moreover, if $T$ is normal, we obtain $[T, (T^*T+sI)^{-1}T^*]=0$ for every $s>0$. Taking $s\to 0$, we get $[T,T^\dag]=0.$
This motivates the following result on normal operators which extends the well-known result on matrices given in  \cite{meyer}:
\begin{prop}\label{normalsdep}
Let $T\in\mathcal{R}(T)$. Then $T$ is  normal operator if and only if $T\in \text{(EP)}\cap \text{(SD)}$.
\end{prop}
\begin{proof}
    For  normal operator $T$  with closed range, we have $T= 0 \oplus T_1$, where $T_1=T_{|N(T)^\perp}$ is the invertible operator on the reducing subspace $N(T)^\perp$. It follows that $T^\dag=0\oplus T_1^{-1}$ and then $T \in \text{(EP)}\cap \text{(SD)}$. Conversely, if $T \in \text{(EP)}\cap \text{(SD)}$, using the identities $T^\dag TT^*=T^* TT^\dag=T^*$  we have
    \begin{eqnarray*}
    T^*T&=&(T^\dag TT^*)T\\
    &=&TT^\dag T^*T\\
    &=&TT^*T^\dag T\\
    &=& TT^*TT^\dag \\
    &=& TT^*.
    \end{eqnarray*}
\end{proof}
 For quasi-normal operators, we have the following result:
 \begin{prop}\label{sdep}
Every quasi-normal operator with closed range is SD.
\end{prop}
\begin{proof}
If $T$ is quasi-normal, then $T$ and $T^*T$ commute and so are $T$ and  $T^*T+sI$. Which yields to $T^*$ and $(T^*T+sI)$ commute for every $s$. Through the limit, we obtain  $T^*$ and $T^\dag$ commute and hence $T$ is $SD.$
\end{proof}
The power of an $(SD)$ is not necessarily $(SD)$. However, since for every quasi-normal operator $T$ and $n\ge 2$,  $T^n$ is quasi-normal, \cite [Proposition 2.6] {emz1}, we derive the following result on powers.
\begin{corol}
  If $T$ is quasi-normal with closed range then $T^{n}$ is SD for every integer $n\ge 0.$
\end{corol}

\begin{Rem}\label{power1}
The general,  reverse order rule for  Moore-Penroses inverse $(AB)^\dag=B^\dag A^\dag$ has been considered by various authors.
  From \cite[Theorem 9]{lu}, if $S, T\in\mathcal{L(H)}$ have closed ranges, then the following statements are equivalent:
\begin{itemize}
 \item[(a)]\;$(ST)^\dag=T^\dag S^\dag$;
 \item[(b)]\;$T(N(ST)^\perp)\subset N(S)^\perp$ and $S^*(N(T^*S^*)^\perp)\subset  N(T^*)^\perp$;
 \item[(c)]\; $S^*S(N(T^*))\subset N(T^*)$ and $TT^*(N(S))\subset N(S)$.
\end{itemize}
It follows that

 $\bullet$  If $U$ is unitary, then $(U^*TU)^\dag =U^*T^\dag U,$

  $\bullet $ $   \omega(ST)=\omega(S)\omega(T) \iff \left\{\begin{array}{l}
       S^*S(N(T^*))\subset N(T^*) \\
   TT^*(N(S))\subset N(S).
    \end{array}\right.$\\

  $\bullet$  $\omega(T^2)=\omega(T)^{2}\iff T^*T(N(T^*))\subset N(T^*)$ and $TT^*(N(T))\subset N(T)$.\\

For an arbitrary operator, the equality $\omega(T^n)=\omega(T)^{n}$ for  $n\in \mathbb{N}$   may fail even for left invertible operators as shown by the  operator $T=S+\frac{1}{2}S^2$ given in \cite [Example 4]{emz2}, where $S$ is   the usual shift on the Hardy space. On the other hand it is clear that when $T$ is EP, then  $\omega(T^n)=\omega(T)^{n}$ for  $n\in \mathbb{N}$.
\end{Rem}
Using the equivalence in Remark \ref{power1}, we derive   the following result on powers.
\begin{prop}
Let $T\in\mathcal{R(H)}$ and $n\geq 2$ be an integer. Suppose that  $T^*T(N(T^{*i}))\subset N(T^{*i})$ and $TT^*(N(T^i))\subset N(T^i)$ for every $i=1,\ldots,n-1$, then $\omega(T^i)=\omega(T)^{i}$ for every $i=1,\ldots,n$.
\end{prop}
\begin{proof}
From \cite[Lemma 3.1]{guy}, for every nonnegative integer $i\le n-1$, we have
$$T^*T(N(T^{*i}))\subset N(T^{*i})\iff T^{*i}T^i(N(T^{*}))\subset N(T^{*}).$$  In this case, $T^2, T^3,\ldots, T^n$ are closed range.   Using  Remark   (\ref{power1}) with $i=1$, gives $(T^2)^{\dag}=(T^{\dag})^2$. Now, for $i=2$ and  $S=T^2$, we use  Remark \ref{power1} again to obtain $T^*T(N(S))\subset N(S)$ and $SS^*(N(T))\subset N(T)$. It follows that   $(TS)^{\dag}=T^{\dag}S^{\dag}$, that is, $T^{3\dag}=T^{\dag}T^{2\dag}=T^{\dag}T^{\dag 2}=T^{\dag 3}$. We  argue inductively to get  $T^{n\dag}=T^{\dag n}$. The result is derived by passing to adjoint.
\end{proof}
\subsection{Basic properties of $n$-EP operators}
The class  of  $n$-\emph{normal operators} for some integer $n\geq 1$ extends the class of normal operators and was    introduced in \cite{jibril}. It  consists of operators $T$ such that   $T^nT^*=T^*T^n$.
In the light of this concept, the class of EP operators is extended  to the class of $n$-EP operators. See  \cite{SLN, wangchun}, for example.

\begin{definition}
For $n\ge 1$, an operator $T\in\mathcal{R(H)}$  is an $n$-EP operator, or  $T\in (n$-EP) if
\begin{equation}\label{n-EP}
    T^nT^\dag=T^\dag T^n.
\end{equation}
\end{definition}
The following examples  show that the class ($n$-EP)  is larger than the
class (EP).
 \begin{example} Let $\mathcal{H}$ be a Hilbert space and  $x, y$ be orthogonal unit vectors. The  partial isometry defined by $T= x\otimes y$ satisfies
 $T^\dag= T^*= y\otimes x$. It follows that
 $TT^\dag=P_{{\mathbb C}x}$ and  $T^\dag T=P_{{\mathbb C}y}$, where $P_{\mathbb{C}x}$ and $P_{\mathbb{C}y}$ are  the orthogonal projections on  $\mathbb{C}x$ and $\mathbb{C}y$ respectively. Thus, $T$ is not EP. Now $T^2=0$ and then  $T$ is $n$-EP for every $n\geq 2$.
\end{example}

 \begin{example}
In the infinite dimension, let $\mathcal{H}=\ell^2$ and $(e_1,e_2,\ldots,e_n,\ldots)$ be standard orthogonal basis for $\ell^2$. Let $\alpha$ be a nonzero real number and define $S$ on $\mathcal{H}$ by
$$Se_k=\begin{cases}
    \alpha e_1,&\quad k=1,\\
     \alpha e_{k+1},&\quad k=2p,\\
    0,&\quad k=2p+1,
\end{cases}$$
where $p=1,2,\ldots.$
Then $S^n=\alpha^n P_{\mathbb{C}e_1}\quad (n\geq 2)$. Simple calculations show that
$S^\dag e_k=\begin{cases}
    \frac{1}{\alpha} e_1,&\quad k=1,\\
     \frac{1}{\alpha} e_{k-1},&\quad k=2p+1,\\
    0,&\quad k=2p.
\end{cases}$
and
$$S^\dagger=\frac{1}{\alpha^2} S^*,\: R(S)=\bigvee_{k\geq 1}\{e_1,e_{2k+1}\}\neq \bigvee_{k\geq 1}\{e_1,e_{2k}\}=R(S^*),$$
and $$S^\dagger S^n=S^nS^\dagger= \alpha^{n-1}P_{\mathbb{C}e_1}\:\:(n\geq 2).$$

So, $S$ is $n$-EP but is not an EP operator.
\end{example}
 It is clear from the definition that  $(1-\text{(EP)}) =\text{(EP)}  \subset(n\mbox{-EP})$ for every $n\geq 2$ and that   $T$ is $n$-EP if and only if $T^*$ is $n$-EP.\\

   Our first observation concerns restrictions to invariant subspace. As for normal and $n-$ normal operator, the restriction to arbitrary invariant subspaces,  the restriction does not preserve this property. On the other hand, the class of normal operators and  $n-$normal operators are closed under unitary equivalence and restrictions to  reducing subspaces. For  $n$-EP class we have the following   properties.
\begin{prop}\label{restriction nEP}Let $T\in\mathcal{R(H)}$ be $n$-EP. Then
    \begin{enumerate}
    \item If  $M$ is  a reducing subspace, then   $T_{|M}$ is $n$-EP.
    \item If $S$ and $T$ are  unitarily equivalent , then  $S\in (n$-EP).
   \end{enumerate}
\end{prop}
\begin{proof}
    (1) We use the identity  $  T^\dag = \lim\limits_{s\to 0}(T^*T+sI)^{-1}T^*$  to show that any reducing subspace is invariant by $  T^\dag$, from which the result is immediate.\\
     (2)  If $S=U^*TU$, we get   $S^\dag=U^*T^\dag U$ and $S^n=U^*T^nU$  and the required result follows.
 \end{proof}
  It is easily seen that a partial isometry $T$ ($T^\dag=T^*$) is $n$-EP if and only it is $n$-normal. The following example shows that a  $n$-EP operator need not be $n$-normal.

\begin{example}
  Consider on ${\mathbb C}^2$  endowed with an orthonormal basis the matrix
  $S=\begin{pmatrix}
2 & 1\\
0 & -2
\end{pmatrix}.$ Since $S^2 = \begin{pmatrix}
4 & 0\\
0 & 4
\end{pmatrix} $,  then  $S^\dag=S^{-1}=\frac{1}{4} S$ and hence $S$ is an EP operator. On the other hand,  it is easily seen that $S$ is not normal and it is $n$-EP for every $n\ge 1$. In particular, $S$ is $2n+1$-EP and is not $2n+1$-power normal for every $n\ge 1$.
  \end{example}

 Before exposing the main properties of $n$-EP operators,  we  show the following result which will be useful later:
\begin{lem}\label{R(An)=R(AnAdag)}
Let  $A\in\mathcal{L(H)}$ be with closed range. For every $n\ge 1,$ we have $$R(A^nA^\dag)=R(A^n) =R(A^nA^*).$$
\end{lem}
\begin{proof}
For the first identity, we write
$$ R(A^n A^\dag)\subset R(A^n)=R(A^nA^\dag A)\subset R(A^nA^\dag).$$
For the second equality, we use the identity $y=A^* u+ v$ with $u\in \mathcal{H}$ and $v\in N(A)$, for an arbitrary $y\in \mathcal{H}$.
\end{proof}
It is known that $T$ is EP if and only if $R(T)=R(T^*)$. For  $(n-\text{EP})$ class, according to \cite[Corollary 2.6]{Mosic}, $T\in(n-\text{EP})$ if and only if $R(T^{n})\subset R(T^*)$ and $R(T^{*n})\subset R(T)$ . In the following, we give a direct and simple proof for this result:
\begin{theo}\label{equivEP}
If $T\in\mathcal{R(H)}$, then $T\in(n-\text{EP})$ if and only if $R(T^{n})\subset R(T^*)$ and $R(T^{*n})\subset R(T)$.
\end{theo}
\begin{proof}
Recall first that $T^{\dag *}=T^{*\dag}.$ Suppose that $T \in(n-\text{EP})$. Using Lemma \ref{R(An)=R(AnAdag)} we get $$R(T^n)=R(T^nT^\dag)=R(T^\dag T^n)\subset R(T^\dag)=R(T^*).$$
Since $ T^{*}\in (n-\text{EP})$, we also have $R(T^{*n})\subset R(T)$.\\
Conversely, assume that $R(T^{n})\subset R(T^*)$ and $R(T^{*n})\subset R(T)$. We have
$$
T^nT^\dag=P_{R(T^*)}T^nT^\dag\\
         =T^\dag T T^nT^\dag\\
         =T^\dag T^n TT^\dag.
$$
Morover, since $T^{*n}=P_{R(T)}T^{*n}=TT^\dag T^{*n}$ we obtain
$T^{n}=T^n TT^\dag$. Therefore, $T^nT^\dag=T^\dag T^n$.
\end{proof}
\begin{corol}
For evry $n\ge 1$, we have  $(n-\text{EP})\subset((n+1)-\text{EP})$.
\end{corol}
\begin{proof}
The result follows from $R(T^{n+1})\subset R(T^n)$ for every $n\geq 0$.
\end{proof}
It has been shown in \cite [Proposition 2.1]{jibril} that $T$ is $n$-normal if and only if  $T^n$ is normal. Concerning $(n-\text{EP})$ classes, we have the following result.
\begin{prop}\label{T^nisEP}
Let $T\in \mathcal{R(H)}$   and $n\geq 1$ such that $T^n $ has closed range.
\begin{enumerate}
\item If $T^n\in$ (EP) then $T\in(n-\text{EP})$.
\item If $T\in$(SD) then $T\in(n-\text{EP})$ if and only if $T^n\in$ (EP).
\end{enumerate}
 \end{prop}
\begin{proof}
$(1)$ Assume that  $T^n\in$ (EP).  Then  $R(T^n)=R(T^{*n})$. It follows that   $R(T^n) \subset R(T^*)$ and  $R(T^{*n})\subset R(T)$ and then $T\in(n-\text{EP})$.\\
$(2)$ Suppose that $T\in(n-\text{EP})\cap \text{(SD)}$ and let  $A=T^n$, $B=T^\dag$.
We have $$ A^*A N(B^*)=T^{*n}T^n N(T^{\dag *})=T^{*n}T^n N(T)=\{0\}\subseteq N(B^{*}).$$
On the other hand, since $T\in(SD)$, for $x\in N(A)$, we obtain
$$ ABB^*x=T^nT^{\dag}T^{\dag *}x=T^\dag T^n T^{\dag *}x=T^\dag T^{\dag *}T^n x=0.$$ Thus $BB^*N(A)\subseteq N(A)$. Now using  Remark \ref{power1}, we derive that $$(T^n T^\dag)^\dag= T^{\dag\dag}T^{n\dag}=TT^{n\dag}.$$ In a similar way, with  $A=T^\dag$ and $B=T^n$, we obtain  $(T^\dag T^n )^\dag=T^{n\dag} T$. Since by assumption $T$ is $n$-EP, we get $TT^{n\dag}=T^{n\dag} T$. This shows that $T^nT^{n\dag}=T^{n\dag} T^n$ or equivalently $T^n$ is $EP.$
\end{proof}
Notice that  if  $T\in$(EP),
$$T^\dag T=(T^\dag T)^2=T^\dag TT^\dag T=T^\dag T^2T^\dag=T(T^\dag)^2T= TT^\dag.$$ From what, we get
 $$T\in(EP) \iff TT^{\dag 2}T=T^\dag T \mbox{ and } T^\dag T^{2}T^\dag=TT^\dag \iff T, T^\dag\in (EP).$$
 For the class ($n$-EP),  starting with  Theorem \ref{equivEP} and  Lemma \ref{R(An)=R(AnAdag)}, we easily arrive at the following result which involve characterizations of $n$-EP.  See also \cite{Mosic}:
\begin{prop}\label{equivEP2}
If $T\in\mathcal{R(H)}$,  then the following assertions are equivalent:
\begin{itemize}
\item[(1)]\;$T\in(n-\text{EP})$.
\item[(2)]\;$T^n=T^\dag T^{n+1}=T^{n+1}T^\dag$.
\item[(3)]\; $TT^{\dag }(T^\dag T^n)^*=(T^\dag T^n)^*$ and $T^\dag T(T^{n}T^\dag)=T^{n}T^\dag$.
\end{itemize}
\end{prop}
Recall that the ascent and the descent of an operator are defined by
$$asc(T)= min \{n\in \mathbb{N} \mbox{ such that } N(T^n)=N(T^{n+1})\}$$ and $$ dsc(T)= min \{n\in\mathbb{N} \mbox{ such that } N(T^n)=N(T^{n+1})\}.$$
An operator has finite ascent ( resp. finite descent) if $asc(T)$ ( resp. $dsc(T)$ is finite). For more details and information on operators with finite ascent and descent, see \cite{Aiena1, Adi} for example.
 \begin{corol}
Let  $A$ be in (EP). For every $n\ge 1,$ we have $$R(A^nA^\dag)=R(A^n) =R(A^{n-1}).$$
\end{corol}
\begin{proof}
For the first identity, we write
$$ R(A^n A^\dag)\subset R(A^n)=R(A^nA^\dag A)\subset R(A^nA^\dag).$$
For the second equality, we use the identities $$ A^nA^\dag=A^{n-1}AA^\dag=A^{n-1}A^\dag A=A^{n-1}.$$
\end{proof}
We also have
\begin{corol} Let  $A$ be $n-$EP. For every $k\ge 0,$ we have $$R(A^{n+k}A^\dag)=R(A^n) =R(A^{n+k}).$$
\end{corol}  Since $T$ is $n$-EP $\iff T^*$ is $n-$EP, and $asc(T^*)=dsc(T)$ the same observations hold for the ascent.
\begin{corol}
Every  $T\in (n$-EP) has finite ascent and descent. Moreover  $asc(T)\le n$ and $dsc(T)\le n.$
\end{corol}
\begin{proof}
     Follows from the identities $T^n=T^{n+1}T^\dag$, $(T^*)^{n}=(T^*)^{n+1}(T^*)^\dag$ and  $asc(T)= dsc(T^*).$
\end{proof}
Recall that an operator $T$ is said to be regular if it has closed range and $N(T)\subset R(T^k)$ for every $k\ge 0.$ For more information about regular operators, see for instance \cite{emz2}. We derive the next result:
\begin{corol}
Let $T$ be in ($n$-EP). If $T$ is regular, then it is invertible.
\end{corol}
\begin{proof}
    By assumption, we have $N(T)\subset R(T^k)$ for every $k\ge 0.$ Now, since $T$ is  $ (n-EP)$, according to Theorem \ref{equivEP} we deduce that $N(T)=\{0\}$ and then $T$ is left invertible. The use of $T^*$  allows to conclude.
\end{proof}
For $n$- normal operators with closed range, we have:
\begin{prop}\label{nN1}
Closed range  $n$- normal operators are $n$-EP.
\end{prop}
\begin{proof}
Suppose that $T$ is $n$-normal.  By Lemma \ref{R(An)=R(AnAdag)},  $$R(T^n)=R(T^nT^*)=R(T^*T^n)\subset R(T^*)$$  and  $$R(T^{*n})=R(T^{*n}T)=R(TT^{*n})\subset R(T).$$ We conclude by appealing Theorem \ref{equivEP}.
\end{proof}
From Proposition \ref{normalsdep},  an operator with closed range is  normal  if and only if it is in the class $\text(EP)\cap \text{(SD)}$. This result is extended  as follows:
\begin{prop}\label{equality1}
If   $T \in$ (SD), then for every  $n\geq 1$,  we have $$ T\in (n-\text{EP})  \iff T  \mbox{ is } n\mbox{-normal}.$$
\end{prop}
\begin{proof}
According to Proposition \ref{nN1}, it suffices to show that  if $T\in  (n-EP)$, then $T $ is   $n$-normal.
Indeed, if  $T\in (n-\mbox{EP})$, by Proposition \ref{equivEP2},  we have $T^{*n}=(T^nTT^\dag)^{*}=TT^\dag T^{*n}$ and since $T \in (SD)$, we obtain,
$$T^{*n}T=TT^\dag T^{*n} T=TT^{*n}T^\dag T=T(T^\dag T T^n)^{*}=TT^{*n}$$
 and the desired result follows.
\end{proof}
\begin{remark}
It follows from  Propositions \ref{T^nisEP} and    \ref{equality1}, that if  $T\in (\text{SD})$, then
$$T\in  (n-\mbox{EP}) \iff   T^n\in(\mbox{EP}) \iff T \mbox{ is } n-\text{normal}.$$
\end{remark}
We finish this part with the following result which generalizes Proposition 5 given in \cite{hartwig} for EP operators:
\begin{prop}
    If $T$ is $n$-EP, then
    $$T^{2n}\ge 0\iff T^{2n-1}\omega(T)\ge 0.$$
\end{prop}
First, note that if $T$ is $n$-EP, then by Theorem \ref{equivEP}, $T^\dag T^k= T^{k-1}=T^kT^\dag$ for all $k\ge n$. In particular, we have $T^\dag T^{2n}= T^{2n-1}=T^{2n}T^\dag$.\\
Assume that $T^{2n}=BB^*$, then $$T^{2n-1}\omega(T)=T^\dag T^{2n}\omega(T)=T^\dag T^{2n}T^{\dag *}=(T^\dag B)(T^\dag B)^*\ge 0.$$
Conversely, assume that $T^{2n-1}\omega(T)= CC^*$.  Then,
$$T^{2n}=T T^{2n-1}TT^\dag=T(T^{2n-1}\omega(T))T^*=TCC^*T^*\ge 0.$$

\section{$n$-hypo-EP operators}
\subsection{Some properties of hypo-EP operators}
The class of hypo-EP operators, (HEP) for short, is introduced by Itoh in \cite{itoh} as an extension    of the EP class. More precisely, this class is defined as the family of operators with generalized inverse $T^\dag$ such that $[T^\dag, T]\ge 0$. The class (HEP)   extends the class of EP operators and satisfies various interesting properties. In particular, from  \cite[Theorem 2.4]{itoh} and \cite{ding1}, for a closed range operators $T$, we have  $$T\in(HEP)\iff T=T^\dag T^2 \iff R(T)\subset R(T^*).$$
The next proposition provides  further characterizations of (HEP) operators:
\begin{prop}
Let $T\in\mathcal{R(H)}$ be a nonzero operator. The following statements are equivalent:
\begin{itemize}
\item[(1)]\;$T\in$ (HEP);
\item[(2)]\;$T^\dag=T^{\dag 2}T$;
\item[(3)]\;$T^\dag T=T^{\dag n}T^n$ for every $n\geq 1$;
\item[(4)]\;  $\|T^{*}x\|\leq c \|Tx\|$ for all $x\in \mathcal{H}$ and for some real number $c\ge 0$.
\end{itemize}
\end{prop}
\begin{proof}
$(1)\iff (2)$ We have
\begin{eqnarray*}
T\in(HEP)&\iff & R(\omega(T))=R(T)\subset R(T^*)\\
         &\iff &\omega(T)=P_{R(T^*)}\omega(T)\\
         &\iff & \omega(T)=T^\dag T\omega(T)\\
         &\iff & T^\dag=(\omega(T))^*=T^\dag (T^\dag T)^*\\
         &\iff& T^\dag= T^{\dag 2}T.
\end{eqnarray*}
$(2)\Rightarrow (3),$ is obvious.\\
$(3)\Rightarrow (1),$ suppose that $T^\dag T=T^{\dag n}T^n$ for every $n\geq 2$. In particular, we have
$T^\dag T=T^{\dag 2}T^2$. Then
\begin{eqnarray*}
T^*T&=&T^*(TT^\dag T)\\
     &=&T^*T(T^{\dag 2} T^2)\\
     &=&(T^*TT^\dag)T^\dag T^2\\
     &=&T^*T^\dag T^2.
\end{eqnarray*}
  Now, let $x\in\mathcal{H}$. We have $$ \|Tx\|^2=\langle T^*Tx,x\rangle=\langle T^*T^\dag T^2x,x\rangle=\langle T^\dag T(Tx),Tx\rangle=\|T^\dag T(Tx)\|^2,  $$ and since $T^\dag T$ is an orthogonal projection, we derive that $T=T^\dag T^2$ and hence $T\in$(HEP).
$(1)\iff (4)$, derives from  the next Douglas factorization theorem
in  \cite{douglas}   $$R(T)\subset R(T^*)\iff \exists c\ge 0 : \, \|T^{*}x\|\leq c \|Tx\| \mbox{ for all } x\in \mathcal{H}.$$
\end{proof}
In the line of \cite[Theorem 5]{camb}, we have:
\begin{theo}
  Let $T\in$ (HEP). Then $$[TT^{\dagger},T+T^{*}]=0 \iff T\in \text{(EP)}.$$
\end{theo}
\begin{proof}
 $\Longrightarrow$:
   We have
  \begin{eqnarray*}
   0=[TT^{\dagger},T+T^{*}]&=& TT^{\dagger}(T+T^{*})-(T+T^{*})TT^{\dagger}\\
    &=&TT^{\dagger}T+TT^{\dagger}T^{*}-TTT^{\dagger}-T^{*}TT^{\dagger}\\
     &=&T+TT^{\dagger}T^{*}-T^{2}T^{\dagger}-T^{*}.
 \end{eqnarray*}
 Multiplying left  by $T^{\dagger}$,  we obtain
 $T^{\dagger}T=(T^{\dagger}T)(TT^{\dagger})$ and by passing to adjoints, $(TT^{\dagger})(T^{\dagger}T)=T^{\dagger}T$. Hence  $P_{R(T)}R(T^*)=R(T^*)$. Therefore $R(T^*)\subset R(T)$.
 Since by assumption we have $R(T)\subset R(T^*)$, we deduce that  $R(T)=R(T^*)$ and therefore $T$ is EP.\\
$\Longleftarrow$: Assume that $T$ is EP. We have $$[TT^{\dagger},T+T^{*}]=T+TT^{\dagger}T^{*}-T^{2}T^{\dagger}-T^{*}=T+T^*-T-T^*=0.$$
\end{proof}
\begin{theo}
   Let $T\in \mathcal{R(H)}$ be such that $[T^{\dagger}T,T+T^{*}]=0$. Then $T$ is hypo-EP.
\end{theo}
\begin{proof}
  We have
  \begin{eqnarray*}
   0&=&[T^{\dagger}T,T+T^{*}]\\
   &=&T^{\dagger}T^{2}+T^{\dagger}TT^{*}-TT^{\dagger}T-T^{*}T^{\dagger}T\\
     &=&T^{\dagger}T^{2}+T^{*}-T-T^{*}T^{\dagger}T .
 \end{eqnarray*}
 Multiplying right by $T^{\dagger}$, we get
 \begin{eqnarray*}
     0= T^{\dagger}T^{2}T^{\dagger}-TT^{\dagger}
 \end{eqnarray*}
and hence $TT^{\dagger}=(T^{\dagger}T)(TT^{\dagger})$ which implies that $R(T)\subset R(T^*)$  and therefore $T$ is $(HEP)$.
\end{proof}
%====================================n-HEP====================================
\subsection{Characterization of $n$-hypo-EP operators}
Taking in count that an operator $T$ with closed range is HEP if and only if $R(T)\subset R(T^*)$,  we extend this notion in the next definition.
\begin{definition}
A closed range operator $T \in\mathcal{L(H)}$  is said to be an $n$-hypo-EP operator ( $n$-HEP, for short) for some  $n\ge 1$,  if $R(T^{n})\subset R(T^*)$.
\end{definition}
It is clear that for every $n\ge 1$, $T$ and $T^*$ are simultaneously in  ($n$-HEP) if and only if $T$ is in ($n$-EP). Moreover, it is easily seen that every left invertible operator is  $n$-HEP. We note in passing that if $T$ is an $n$-HEP regular operator then it is left invertible.
 \begin{remark} The class of $n$-HEP operators on $\mathcal{H}$,   noted by ($n$-HEP),  is closed under unitary equivalence and restriction on reducing subspaces. Moreover,
 \begin{itemize}
 \item[(1)]\; (1-HEP)=(HEP);
 \item[(2)]\;$(n-\mbox{EP})\subset (n-\mbox{HEP})$;
 \item[(3)]\;$(n-\mbox{HEP})\subset((n+1)-\mbox{HEP})$, for every $n\ge 1$.
 \end{itemize}
 \end{remark}

As mentioned before,  $T\in (\mbox{HEP}) \iff T=T^\dag T^{2}\iff T^*=T^*T^\dag T$.
We extend next  the previous characterization to $(n-HEP)$. The proof run in a similar way and is omitted.
\begin{prop}\label{caractnHEP}
For  $T \in\mathcal{R(H)}$, the following  are equivalent:
\begin{itemize}
\item[(1)]\;$T\in$ ($n$-HEP);
\item[(2)]\;$T^n=T^\dag T^{n+1}$;
\item[(3)]\;$T^n T^\dag=T^\dag T^{n+1} T^\dag$;
\item[(4)]\; $T^{*n}=T^{*n}T^\dag T$;
\item[(5)]\;$\|T^{*n}x\|\leq c \|Tx\|$ for some $c\geq 0$ and all $x\in\mathcal{H}$.
\item[(6)]\;\;$\|T^{*n}x\|\leq c \|\omega(T)x\|$ for some $c\geq 0$ and all $x\in\mathcal{H}$.
\end{itemize}
\end{prop}
 We also have:
\begin{prop}
For  $T \in\mathcal{R(H)}$, if $T^n\in (HEP)$ then $T\in (n-\mbox{HEP})$.
 \end{prop}
\begin{proof}
Suppose that $T^n\in (HEP)$, then $R(T^n)\subset R(T^{*n})\subset R(T^*)$.
\end{proof}

\begin{prop}
Let $T\in\mathcal{R}(H)$ be such that one of the following conditions is satisfied:
\begin{itemize}
    \item [i.] $[T^\dag T,T^n+T^\dag]=0$;
    \item [ii.] $[T^\dag T,T^n+T^*]=0$,
\end{itemize}
then $T\in(n-HEP)$.
\end{prop}
\begin{proof}
 Suppose that (i) is satisfied. Then,
\begin{eqnarray*}
0&=&[T^\dag T,T^n+T^\dag]\\
 &=&T^\dag T^{n+1}+T^\dag TT^\dag-T^nT^\dag T-T^{\dag2}T\\
 &=&T^\dag T^{n+1}+T^\dag-T^n-T^{\dag2}T
\end{eqnarray*}
Multiplying right by $T^\dag$, we get $$T^\dag T^{n+1}T^\dag=T^nT^\dag.$$
By Proposition \ref{caractnHEP} we conclude that $T\in(n-\mbox{HEP})$.\\
If (ii) is true, then
\begin{eqnarray*}
0&=&[T^\dag T,T^n+T^{*}]\\
 &=&T^\dag T^{n+1}+T^\dag TT^{*}-T^nT^\dag T-T^{*}T^\dag T\\
 &=&T^\dag T^{n+1}+T^{*}-T^nT^\dag T-T^{*}T^{\dag}T
\end{eqnarray*}
Multiplying right by $T^\dag$, we have $$T^\dag T^{n+1}T^\dag=T^nT^\dag.$$
By Proposition \ref{caractnHEP}, we derive that $T\in(n-HEP)$.
\end{proof}
\begin{prop}
Let $T\in\mathcal{R(H)}$  be $n$-HEP. If $T$ satisfies one of the following statements,
\begin{itemize}
    \item [i.] $[TT^\dag ,T^n+T^\dag]=0$;
     \item [ii.] $[TT^\dag ,T^n+T^{*}]=0$;
     \item [iii.]$[T ,T^nT^\dag]=0$,
\end{itemize}
then $T\in(n-\mbox{EP})$.
\end{prop}
\begin{proof} (i): It suffices to show that $T^*\in(n-\mbox{HEP})$.
 We have,
$$ 0= [TT^\dag, T^n+T^\dag]= T^n +TT^{\dag 2}- T^{n+1}T^\dag-T^\dag.$$
Multiplying right by $T$, we get $$(TT^{\dag })(T^\dag T)=T^\dag T.$$
Thus, $R(T^*)=R(T^\dag T)\subset R(T)$ and hence $T^*\in(HEP)\subset(n-HEP)$.\\
 (ii): As before, $$ 0= [TT^\dag,T^n+T^{*}]=T^n +TT^\dag T^{*}-T^{n+1}T^\dag-T^{*}.$$
 Multiplying left by $TT^\dag$, we get $ T^n = T^{n+1}T^\dag$ and  by  Proposition \ref{equivEP2}, we conclude that $T$ is $n$-EP.\\
(iii): Again we have,
$0= [T ,T^nT^\dag]=T^{n+1}T^\dag-T^n$ and by
Proposition \ref{equivEP2}, we derive that $T$ is $n$-EP.
\end{proof}

 \end{document}